\documentclass[12pt]{article}
\def\ni{\noindent}
\def\i{\indent}
\def\l{\ell}
\def\t{\theta}
\begin{document}
\begin{center}
{\bf MULTIVARIATE PREDICTION AND MATRIX SZEG\"O THEORY}
\end{center}
\begin{center}
{\bf N. H. BINGHAM}
\end{center}
\begin{center}
{\bf Abstract}
\end{center}
\ \ \ \ \  Following the recent survey by the same author of Szeg\"o's theorem and orthogonal polynomials on the unit circle (OPUC) in the scalar case, we survey the corresponding multivariate prediction theory and matrix OPUC (MOPUC).\\

\ni {\it AMS 2000 subject classifications}.  Primary 60G10, secondary 60G25 \\

\ni {\it Key words and phrases} Stationary process, vector-valued, multivariate prediction theory, multivariate orthogonal polynomials on the unit circle (MOPUC), Kolmogorov isomorphism theorem, Verblunsky's theorem, Szeg\"o's theorem \\

\begin{center}
CONTENTS
\end{center}
\ni \S 0.  Introduction \\
\ni \S 1.  The Kolmogorov Isomorphism Theorem \\
\ni \S 2.  Verblunsky's theorem \\
\ni \S 3.  Szeg\"o's theorem \\
\ni \S 4.  Matrix spectral factorisations and matrix Szeg\"o functions \\
\ni \S 5.  The strong Szeg\"o theorem \\
\ni \S 6.  Baxter's inequality and Baxter's theorem \\
\ni \S 7.  Nehari sequences and the Levinson-McKean condition \\
\ni \S 8.  Pure minimality \\
\ni \S 9.  Positive angle and the matrix muckenhoupt condition \\
\ni \S 10. Complete regularity \\
\ni \S 11. Hankel operators \\
\ni \S 12. Open questions \\
References \\

\newpage

\ni {\bf \S 0.  Introduction} \\

\i Prediction theory for discrete-time stationary stochastic processes, a topic in time series, has been long studied and has benefited greatly from
the recent work on orthogonal polynomials on the unit circle (OPUC); see e.g. Simon's books [Si1], [Si2], [Si3] for background and references.
Here the partial autocorrelation function (PACF) plays an important role, as it provides an unrestricted parametrization of the relevant spectral
measure $\mu$; see e.g. the papers by Inoue [In1], [In2], [In3], Inoue and Kasahara [InKa1], [InKa2] and the survey by Bingham [Bi] for background and references.
The PACF is essentially the sequence of Verblunsky coefficients, and the unrestricted parametrization is Verblunsky's theorem, in the language of OPUC; see [Si1], Ch.
1,2.\\
\i While this is very satisfactory for univariate time series, one often encounters multivariate time series, particularly in areas such as mathematical finance, where the dimensionality
$\l$ is the number of risky assets held in a portfolio; $\l$ reflects the need to diversify one's portfolio (in the style of Markowitz), and may be large (see e.g. [BiFrKi] for a case in point).  Multivariate time series form an important
area (see e.g. the books by Hannan [Ha] and Reinsel [Re]).  Likewise, multivariate prediction goes back to work by Wiener and Masani [WiMa1], [WiMa2], Helson and Lowdenslager [HeLo] of 1957-61
 (also Masani [Mas1] -- [Mas5]).  Just as in the univariate case with OPUC, in the multivariate case the matrix theory of OPUC -- MOPUC for short below -- is crucially
 relevant.  Following the great stimulus to OPUC provided by Simon's books, MOPUC has recently been extensively developed.  Our purpose here is to survey these developments
with a view to multivariate prediction theory, providing a matrix sequel to [Bi]. \\

\ni {\bf 1.  The Kolmogorov Isomorphism Theorem} \\

\i Stone's theorem ([Sto], [RiNa] \S 137, extended to semigroups in [HiPh] XXII) tells us that a group $U = (U_t)$ of unitary transformations has a spectral representation of the form
$$
U_t = \int e^{i \t t} E(d \t),
$$
where $E(.)$ is a projection-valued random measure.  For a stationary stochastic process $X = (X_t)$ (here time $t$ is discrete and $X_t$ is a random complex $\l$-vector), write $U$ for the shift $t \mapsto t+1$.  Since we can write $X_t = U^t X_0$, this gives us a spectral representation for $X = (X_t)$:
$$
X_t = U^t X_0 = \int e^{i \t t} E(d\t) X_0
$$
(see e.g. [Rozanov [Ro2], Th. I.4.2).  On the left, we have explicit dependence on time $t$ -- we are in the {\it time domain}.  On the right, for fixed $t$ we have dependence on $\t$ -- we are in the {\it frequency domain}.  The spectral representation may be summarized as
$$
X_t \leftrightarrow e^{it.},
$$
which expresses the {\it Kolmogorov Isomorphism Theorem} (Kolmogorov [Ko] in 1941; cf. [Ro2], 17, 33, Masani [Mas4] \S \S 6,7).  Multiplying these expressions for $X_s$, $X_{s+t}$ and using stationarity, we get on taking expectations the spectral representation of the correlation matrix:
$$
{\gamma}_n = \int e^{- i n \t} \mu(d \t),
$$
where $\mu$ is the {\it spectral measure}.  This is {\it Herglotz's theorem}; see e.g. [Ro2], 19-20. \\

\ni {\bf \S 2.  Verblunsky's theorem}\\

\i In the univariate case, the Verblunsky coefficients $a = (a_n)_{n=1}^{\infty}$ satisfy $|a_n| < 1$, and the bijection
$$
a \leftrightarrow \mu
$$
of Verblunsky's theorem provides the unrestricted parametrization so useful in statistics and prediction theory; see e.g. Pourahmadi [Pou1], [Pou2],  [Bi] \S 2 for background here.\\
\i In the $\l$-dimensional case, where one uses MOPUC rather than OPUC, the spectral measure $\mu$ is now $\l \times \l$ matrix-valued, and corresponds to the $\l \times \l$ covariance matrix by Herglotz's theorem.  The Verblunsky coefficients
have been studied by Damanik, Pushnitski and Simon [DaPuSi], \S 3.  They show (3.10) that the Verblunsky coefficients are now $\l \times \l$ matrices on the unit circle, satisfying
$$
\Vert a_n \Vert < 1,
$$
that any sequence $a$  of such matrices can arise in this way, and that the map $a \leftrightarrow \mu$ is again a bijection (Verblunsky's theorem for MOPUC -- [DaPuSi], Th. 3.12).  Their method is Bernstein-Szeg\"o
approximation  ([DaPuSi] \S 3.6; cf. [Si1] Th. 1.7.8 in the scalar case).\\
\i Results of this type go back in the statistical literature to Morf, Vieira and Kailath [MoViKa] in 1978.  They showed that the $a_n$ have singular values ${\alpha}_n$ with $|{\alpha}_n| \leq 1$. \\
\i The Szeg\"o recursion that leads to OPUC is known in the time-series literature as the {\it Levinson-Durbin algorithm}.  The Levinson-Durbin algorithm was extended to the multivariate case by Whittle [Wh] in 1963, and by Wiggins and Robinson [WiRo] in 1965.  In [MoViKa], the authors use the results of [Wh] and [WiRo], which they call the LWR algorithm, and give a normalized version of it.  They show that unless the Levinson-Durbin algorithm stops (corresponding to some correlation matrix failing to be positive-definite), $|{\alpha}_n| < 1$ for all $n$. \\

\ni {\bf \S 3.  Szeg\"o's theorem}\\

\i Derevyagin, Holtz, Khrushchev and Tyaglov [DeHoKhTy] continue this study, again using Bernstein-Szeg\"o approximation (\S \S 5,6).  They show (\S 7) that (using $\dagger$ for the adjoint
(conjugate transpose),
${\mu}'$ for the density of (the absolutely continuous component of) $\mu$, denoted $w$ below)
$$
 \log \Pi_{n=1}^{\infty} det(1 - a_n^{\dagger} a_n) = \int tr \ \log {\mu}' d\theta/2 \pi
 = \int tr \ \log w d\theta/2 \pi
 $$
 for any non-trivial (i.e. of infinite support) matrix-valued probability measure on the unit circle.  Call those for which the integral on the right is $> -\infty$ {\it Szeg\"o measures}, and the finiteness of the integral {\it Szeg\"o's condition}, $(Sz)$.  They deduce that $\mu$ is a Szeg\"o measure iff
 $$
 \sum \Vert a_n^{\dagger} a_n \Vert < \infty
 $$
 ("$a \in L_2$").  This is Szeg\"o's theorem for MOPUC.  They continue (\S 8) with the matrix version of the Helson-Lowdenslager theorem:
 $$
 \exp \int \frac{1}{\l} tr \ \log {\mu}' d\theta / 2 \pi =
{\inf}_{A, P} \int \frac{1}{\l} tr [(A + P)^{\dagger} d \mu (A + P)],
$$
where the infimum is over all matrices $A$ of determinant 1 and all trigonometric polynomials $P(e^{i \theta}) = {\sum}_{k > 0} A_k e^{i k \theta}$. \\
\i The Wold decomposition extends from the scalar to the vector case.  See e.g. Masani [Ma1], [Mas4], Hannan [Ha], and in the context of operator theory, Sz.-Nagy and Foias [SzNF], Ch. 1, Nikolskii [Nik] Ch. 1.  As in the scalar case, $(Sz)$ is the condition that the deterministic component in the Wold decomposition should vanish, leaving only the moving-average component.  This is the condition of {\it non-determinism}, $(ND)$ ("$(Sz) = (ND)$"). \\
\i As in the scalar case, it is often useful to strengthen the Szeg\"o condition by requiring explicitly that the singular component ${\mu}_s$ of the spectral measure $\mu$ should vanish.  This is called the condition of {\it pure non-determinism}, $(PND)$:
$$
(PND) = (ND) + \{{\mu}_s = 0 \} = (Sz) + \{ {\mu}_s = 0 \}.
$$

\ni {\bf \S 4.  Matrix spectral factorizations and matrix Szeg\"o functions} \\

\i Factorizations are already present in the scalar case.  For an analytic function in the Hardy space on the disc, identify the boundary values of the function on the unit circle with the function itself, as usual; then
the spectral density $w$ and the Szeg\"o function $h$ are related by
$$
w = h \bar h = |h|^2.
$$
Here $h$ is in the Hardy space $H_2$, and is an outer function (see e.g. [Bi] for references and details); one can think of $h$ as the `analytic square root' of $w$. \\
\i In the matrix case, one speaks of the {\it spectral factorization problem}.  The two main traditional approaches to multivariate prediction theory are those of Helson and Lowdenslager [HeLo], based on approximation by polynomials, and of Wiener and Masani [WiMas1] -- [WiMas3], based on matrix factorization:
$$
W = G G^{\dagger},           \eqno(WM)
$$
where again $G$ is an outer function in the sense of matrix-valued Hardy spaces (see e.g. [Pel3], \S 13.3 and Appendix 2.3, or [RoRo], Ch. 4,5, for definitions and details; cf. Stein [Ste], Ch. III).\footnote{Wiener and Masani [WiMas3], Def. 3.5 use the term optimal in place of outer.  Rozanov [Ro2], Th. II.4.2 uses the term maximal.}  We refer to factorizations of the form $(WM)$ as {\it Wiener-Masani factorizations}, and to $G$, the matrix analogue of the Szeg\"o function $h$ in the scalar case ($\Delta$ in the notation of [Si1]), as a {\it (matrix) Szeg\"o function}.  Wiener-Masani factorizations are unique up to unitary equivalence ([WiMas2], Th. 8.12). \\
\i This work is well summarized in Masani [Mas3].  Here one finds: \\
(i) the Wold decomposition (\S 4; [WiMas1], Th. 7.11).  In the scalar case, the deterministic (time-independent) component of $X_t$ corresponds to the singular component ${\mu}_s$ of the spectral measure, the moving-average component to the absolutely continuous component ${\mu}_{ac}$ in the Lebesgue decomposition.  In the matrix case, the terms Wold-Zasuhin decomposition and Lebesgue-Cram\'er decomposition are often used; the correspondence is exact in the full-rank case (below), but not in general [Mas1];  \\
(ii) the Kolmogorov Isomorphism Theorem, between the time domain and the spectral domain (\S \S 6, 7); \\
(iii) Wiener-Masani factorizations $(WM)$ (Th. 9.7 -- see also Rozanov [Ro1], [Ro2]); \\
(iv) the matrix extension of Kolmogorov's formula for the one-step prediction error (eq. (10.1), the main result of [WiMas1] (Th. 7.10); see also Whittle [Wh]);  \\
(v) convergence of the finite-past predictor to the infinite-past predictor (\S 13) -- cf. Baxter's inequality, \S 6 below). \\
To these, we also add \\
(vi) Whittle's multivariate extension of the Levinson-Durbin algorithm, mentioned in \S 2. \\
See also [Mas4] for extensive commentary on Wiener's work in this area. \\
\i In the matrix case, one needs to discriminate between the {\it full-rank} and {\it degenerate-rank} cases -- where the rank of the spectral density matrix $W$ is full ($\ell$) or degenerate ($m < \ell$).  The degenerate-rank case is considered in [Ro2], [Mas3], \S \S 11,12, [WiMas3] for $\ell = 2$, Matveev [Mat] in the general case; we refer there for details.  The generic, and easier, case is the full-rank case, where $\Gamma$ is positive-definite (the contrast between the full- and degenerate-rank cases is similar to that arising in, e.g., regression, where one encounters multi-collinearity; see e.g. [BiFr], \S 7.4. \\
\i This interesting work dates from the 1960s, before the work of Fefferman and Stein on BMO and of Sarason on VMO.  Armed with these, Peller [Pel1] in 1990 considered
matrix spectral factorizations
$$
w = h^{\ast} h = h_{\sharp} h_{\sharp}^{\ast}
$$
(recall $h$ is determined to within unitary equivalence).  He introduced the {\it phase function}
$$
u = h_{\sharp}^{\ast} h^{-1},
$$
the analogue of the phase function $\bar h/h$ of \S 5 in the scalar case.  He showed, among other results, that the process is completely regular (\S 10 below) iff
$$
u \in VMO.
$$
\i Arov and Dym [ArDy3], \S 3.16] give matrix factorizations of positive definite functions into factors from the Nevanlinna class. \\
\i The simplest, and principal, case is that of a purely non-deterministic process of full rank.  See e.g. several papers by Ephremidze, Janashia andLagvilada [EpJaLa], [EpLa], [JaLaEp].  In particular, a (matrix) function in a Hardy space is outer iff its (scalar) determinant is outer.k,   \\
\i More general than the matrix-valued case is the operator-valued case.  This is important in operator theory, in non-commutative probability theory, non-commutative Hardy-space theory, non-commutative martingale inequalities etc.  See e.g. [RoRo], esp. Ch. 6, Curtain and Zwart [CuZw], Barclay [Bar1], [Bar2], Mei [Mei], Peller [Pel1] -- [Pel3]. \\

\ni {\bf \S 5.  The strong Szeg\"o theorem} \\

\i The strong Szeg\"o theorem, as presented in e.g. [Si1] Ch. 6, [Bi] \S 5, extends in full to the matrix case.  For a short proof, see B\"ottcher [Bo1]; cf. [Bo2], [Bo3], Basor and Widom [BasWi].
A different approach has been given more recently by Chanzy [Cha1], [Cha2]. \\

\ni {\bf \S 6.  Baxter's inequality and Baxter's theorem} \\

\i  Baxter used OPUC in a series of probabilistic papers of 1961-63 ([Bax1] -- [Bax3]).  His results concern, among other things, the weak and strong forms of Szeg\"o's limit theorem (for Toeplitz determinants), finite and
infinite Wiener-Hopf equations (in discrete time $n = 0,1,2,\ldots$: finite with $\sum_{k=0}^n$, infinite with $\sum_{k=0}^{\infty}$), and the convergence of finite-predictor coefficients (in which one is given a
finite section of the past of length $n$) to the corresponding infinite-predictor coefficients.  This last depends on {\it Baxter's inequality} [Bax3].  Baxter's inequality was used by Simon [Si1], Ch. 5, in his proof that the
Verblunsky coefficients $a \in {\l}_1$ iff the correlation function $\gamma \in {\l}_1$, the spectral measure $\mu$ is absolutely continuous, and its density ${\mu}' = w$ is continuous and positive.  Simon calls this result
{\it Baxter's theorem} (though Baxter did not formulate the result in this form, and `the Baxter-Simon theorem' might be better here, but we will follow [Si1]).  Perhaps because of this rather involved history, and the fact that Simon's book [Si1] is still comparatively recent, there is as yet no matrix extension to Baxter's theorem.  We raise here the question of obtaining one, and turn now to what is known. \\
\i Baxter's inequality and convergence of finite predictors in the matrix case were considered by Masani in 1966 ([Ma3] \S 13) and by Cheng and Pourahmadi [ChPo] in 1993.  Theoretical progress in the area since then has been extensive, and the question arises of weakening the conditions that they impose.  For more recent developments in the
scalar case, see [InKa2]. \\
\i Results on approximation by such finite-section operators have been given in great generality by Seidel and Silbermann [SeSi] (see \S 2.5.4), using Banach-algebra techniques (as did Baxter and Simon). \\

\ni {\bf \S 7.  Nehari sequences and the Levinson--McKean condition} \\

\i Nehari's theorem of 1957 states that a Hankel operator is a bounded map from ${\l}_2$ on the natural numbers to itself iff the sequence generating it is the sequence of negative Fourier coefficients of a bounded function.  See e.g. [Si1] Th. 6.2.17
("The modern theory of Hankel operators started with the following result of Nehari"), Peller [Pel3].  Finding such a generating sequence is thus a type of {\it moment problem}, and as with other moment problems there may be no solution, a unique
 solution or more than one solution; the moment problem is then called insoluble, determinate or indeterminate.  The indeterminate case is particularly important; the generating sequence is then called a {\it Nehari sequence}.  This Nehari moment
(or interpolation) problem was considered by Adamjan, Arov and Krein [AdArKr] in 1968; they described the solution set in terms of Sarason's concept of {\it rigidity} [Sa1].  Rigidity has previously been studied in this area in connection with the concept of
 {\it complete non-determinism} (CND); see Bloomfield, Jewell and Hayashi [BlJeHa].  It turns out that (CND) is equivalent to the {\it intersection of past and future} property (IPF) [IK2].  It has very recently been shown that both are equivalent to the {\it Levinson--McKean property}
(the name comes from work of Levinson and McKean, [LeMcK] p. 105, in continuous time; the name is given by analogy in the discrete-time case).  Phase functions (\S4) play a crucial role here; see Kasahara and Bingham [KaBi] for details. \\
\i The question arises of matrix extensions of these results (work in progress).  The matrix Nehari problem has been considered in detail by Arov and Dym [ArDy1], [ArDy2], [ArDy3] Ch. 4, 7, 10 (`strong regularity').  See also Dym's review of Peller's book on Hankel operators [Pel3] (MR1949210 (2004e:47040)).  For matrix phase functions, see \S 4.  \\
\i Related to this Nehari problem is the Schur (interpolation) problem, the matrix case of which is considered at book length in Dubovoj, Fritzsche and Kirstein [DuFrKi]; cf. [ArDy], \S 7.6. \\

\ni {\bf \S 8.  Pure minimality} \\

\i In the scalar case, pure minimality is characterized by (${\mu}_s = 0$ and) Kolmogorov's condition $1/w  \in L_1$, using $w$ for the density of the spectral measure $\mu$ (now absolutely continuous).  This result extends to the multivariate case;
see Makagon and Weron [MaWe], [Pou3, Th. 8.10].  The spectral density is now a matrix $W$, and its inverse $W^{-1}$ is integrable.  It is thus natural that the condition
$$
W^{-1} \in L_1
$$
should be imposed in studying processes subject to stronger regularity conditions than pure minimality.  We turn below to two such conditions -- positive angle (\S 9) and complete regularity (\S 10).  We regard the four conditions in \S\S 7-10, which are in increasing order of strength, as {\it intermediate} conditions, being intermediate between the {\it weak} conditions (ND), (PND) (Szeg\"o condition + ${\mu}_s = 0$) and the strong conditions (B), (sSz) (Baxter's condition and the strong Szeg\"o condition).  See e.g. [Bi] for the scalar case. \\

\ni {\bf 9.  Positive angle and the matrix Muckenhoupt condition} \\

\i The {\it Muckenhoupt condition} $(A_2)$ of analysis is important in many areas; see e.g. [Bi] \S 6.2 for background and references.  It occurs in connection with the  {\it positive angle} condition, $(PA)$, and the conditions of Helson and Szeg\"o [HeSz] and
of Helson and Sarason [HeSa].  Matrix versions of the Muckenhoupt condition are considered at length by Arov and Dym [ArDy1], [ArDy3].  For matrix versions of the Helson-Szeg\"o condition, see Pourahmadi [Pou1]. \\
\i Treil and Volberg [TrVo1] show that the following matrix Muckenhoupt condition is necessary and sufficient for the positive-angle condition $(PA)$ in the multivariate case:
$$
{\sup}_{I} \Vert \Bigl( \frac{1}{|I|} \int_I W \Bigr)^{1/2}  \Bigl( \frac{1}{|I|} \int_I W^{-1} \Bigr)^{1/2} \Vert < \infty, \eqno(A_2)
$$
where the supremum is taken over all intervals $I$ of the unit circle.  Here the condition that $W$ be invertible a.e. (which corresponds to the `pure' in pure minimality, see \S 8 above) has to be imposed explicitly, as noted by Peller in his review of [TrVo2] (MR1428818 (99k;42073)). \\
\i As shown in [HeSz], [HeSa], the condition $(PA)$ (and so $(A_2)$ by above) is equivalent to a condition on the sequence $\rho(n)$ of regularity coefficients of the form
$$
\rho(.) < 1.
$$
We note that in his review of earlier work on this problem by Makagon, Miamee and Schr\"oder [MaMiSc], Pourahmadi says "Attempts to obtain a similar result for $q$-variate stationary sequences have been unyielding" (MR1443841 (98g:60074)). \\

\ni {\bf 10.  Complete regularity} \\

\i We turn now to a strengthening of the conditions of \S 9 above.  The process is said to be {\it completely regular} if $\rho(n) \to 0$ as $n \to \infty$; see [IbRo], Ch. 4, 5.
 It was shown by Treil and Volberg [TV2] that complete regularity is equivalent to the following strengthening of the Muckenhoupt condition $(A_2)$:
 $$
{\limsup}_{|I| \to 0} \Vert \Bigl( \frac{1}{|I|} \int_I W \Bigr)^{1/2}  \Bigl( \frac{1}{|I|} \int_I W^{-1} \Bigr)^{1/2} \Vert < \infty
$$
(cf. Peller [Pel1], [Pel3]); here as before $W^{-1} \in L_1$ needs to be assumed explicitly.  Note the form ("$\rho(.) \to 0$, $\limsup ... = 1$") of the strengthenings here of the conditions ("$\rho(.) < 1$, $\sup ... < \infty$") of \S 7 above.   \\

\ni {\bf \S 11.  Hankel operators} \\

\i Prediction theory has always involved Toeplitz operators (as in the book [GrSz] by Grenander and Szeg\"o), and Toeplitz and Hankel operators have many links in operator theory.  So it is natural that Hankel operators are useful in prediction theory.  A monograph treatment of Hankel operators is given by Peller ([Pe3]; see also the review by Dym cited above in \S 7).  Connections of Hankel operators with the matrix Muckenhoupt condition and with the matricial Nehari problem are considered by Arov and Dym in [ArDy3], Ch. 10, 11.\\

\newpage

\ni {\bf \S 12.  Open questions} \\

\i We mention two. \\
{\it Question 1}.  Find the matrix version of Baxter's theorem. \\
\i As mentioned in \S 6, the matrix version of Baxter's inequality provides a good starting-point. \\
{\it Question 2}.  Find the matrix version of [KaBi]. \\
\i This hinges on solution of the matrix Nehari problem -- the step
$$
\Gamma \to H.
$$
We hope to return to this elsewhere. \\

\begin{center}
{\bf References}
\end{center}
\ni [AdArKr] V. M. Adamjan, D. Z. Arov and M. G. Krein, Infinite Hankel matrices and generalized problems of Carath\'eodory--Fej\'er and I. Schur.  {\sl Functional Anal. Appl.} {\bf 2} (1968), 269-281. \\
\ni [ArDy1] D. Z. Arov and H. Dym, Matricial Nehari problems, $J$-inner functions and the Muckenhoupt condition.  {\sl J. Funct. Anal.} {\bf 181} (2001), 227-299. \\
\ni [ArDy2] D. Z. Arov and H. Dym, Criteria for the strong regularity of $J$-inner functions and $\gamma$-generating functions.  {\sl J. Math. Anal. Appl.} {\bf 280} (2003), 387-399. \\
\ni [ArDyi3] D. Z. Arov and H. Dym, {\sl $J$-contractive matrix valued functions and related topics}, Encycl. Math. Appl. 116, Cambridge Univ. Press, Cambridge, 2008. \\
\ni [Bar1] S. J. Barclay, Continuity of the spectral factorization maping.  {\sl J. London Math. Soc.} {\bf 70} (2004), 763-779. \\
\ni [Bar2] S. J. Barclay, A solution of the Douglas-Rudin problem for matrix-valued functions.  {\sl Proc. London Math. Soc.} {\bf 99} (2009), 757-786. \\
\ni [BasWi] E. L. Basor and H. Widom, On a Toeplitz determinant identity of Borodin and Okounkov.  {\sl Integral Equations and Operator Theory} {\bf 37} (2000), 397-401.\\
\ni [Bax1] G. Baxter, A convergence equivalence related to polynomials orthogonal on the unit circle.  {\sl Trans. Amer. Math. Soc.} {\bf 99} (1961), 471-487.\\
\ni [Bax2] G. Baxter, An asymptotic result for the finite predictor.  {\sl Math. Scand.} {\bf 10} (1962), 137-144.\\
\ni [Bax3] G. Baxter, A norm inequality for a "finite-section" Wiener-Hopf equation.  {\sl Illinois J. Math.} {\bf 7} (1963), 97-103.\\
\ni [Bi] N. H. Bingham, Szeg\"o's theorem and its probabilistic descendants,\\
arXiv:1108.0368v1 [math.PR] 1 Aug 2011.\\
\ni [BiFr] N. H. Bingham and J. M. Fry, {\sl Regression: Linear models in statistics}.  SUMS (Springer Undergraduate Mathematics Series), 2010. \\
\ni [BiFrKi] N. H. Bingham, J. M. Fry and R. Kiesel, Multivariate elliptic processes.  {\sl Statistica Neerlandica} {\bf 64} (2010), 352-366. \\
\ni [BiInKa] N. H. Bingham, A. Inoue and Y. Kasahara, An explicit representation of Verblunsky coefficients.
{\sl Statistics and Probability Letters}, \\
http://dx.doi.org/10.1016/j.spl.2011.11.004. \\
\ni [BlJeHa]P. Bloomfield, N. P. Jewell and E. Hayashi, Characterization of completely non-deterministic stochastic processes.  {\sl Pacific J. Math.} {\bf 107} (1983), 307-317. \\
\ni [Bo1] A. B\"ottcher, One more proof of the Borodin-Okounkiv formula for Toeplitz determinants.  {\sl Integral Equations and Operator Theory} {\bf 41} (2001), 13 -125.\\
\ni [Bo2] A. B\"ottcher, Featured review of the Borodin-Okounkov and Basor-Widom papers.  {\sl Mathematical Reviews} 1790118/9 (2001g:47042a,b).\\
\ni [Bo3] A. B\"ottcher, On the determinant formulas by Borodin, Okounkov, Baik, Deift and Rains.  {\sl Operator Th. Adv. Appl.} {135}, 91-99, Birkh\"auser, Basel, 2002.\\
\ni [Cha1] J. Chanzy, Th\'eor\`emes-limite de Sz\"ego dans le cas matriciel.  {\sl Proc. Japan Acad. A} {\bf 82} (2006), 113-116. \\
\ni [Cha2] J. Chanzy, Inversion d'un op\'erateur de Toeplitz tronqu\'e \`a symbole matriciel et th\'eor\`emes-limite de Szeg\"o.  {\sl Ann. Math. Blaise Pascal} {\bf 13} (2006), 111-205. \\
\ni [ChePou] R. Cheng and M. Pourahmadi, Baxter's inequality and convergence of predictors of multivariate stochastic processes.  {\sl Probab. Th. Rel. Fields} {\bf 95} (1993), 115-124. \\
\ni [CuZw] R. F. Curtain and H. Zwart, {\sl An introduction to infinite-dimensional linear systems}.  Springer, 1995. \\
\ni [DaPuSi] D. Damanik, A. Pushnitski and B. Simon, The analytic theory of matrix orthogonal polynomials.  {\sl Surveys in Approximation Theory} {\bf 4} (2008), 1-85. \\
\ni [DeHoKhTy] M. Derevyagin, O. Holtz, S. Khrushchev and M. Tyaglov, Szeg\"o's theorem for matrix orthogonal polynomials.  arXiv:1104.4999v1 [math.CA] 26 April 2011. \\
\ni [Dev] A. Devinatz, The factorization of operator-valued functions.  {\sl Ann. Math.} {\bf 73} (1961), 458-495. \\
\ni [DuFrKi] V. K. Dubovoj, B. Fritzsche and B. Kirstein, {\sl Matricial version of the classical Schur problem}.  Teubner, Stuttgart, 1992. \\
\ni [EpJaLa] L. Ephremidze, G. Janashia and E. Lagvilava, An analytic proof of the matrix spectral factorization theorem.  {\sl Georgian Math. J.} {\bf 15} (2008), 241-249. \\
\ni [EpLa] L. Ephremidze and E. Lagvilada, Remark on outer analytic matrix functiions.  {\sl Proc. A. Razmadze Math. Inst.} {\bf 152} (2010), 29-32. \\
\ni [GrSz] U. Grenander and G. Szeg\"o, {\sl Toeplitz forms and their applications}.  Univ. California Press, Berkeley, CA, 1958. \\
\ni [Ha] E. J. Hannan, {\sl Multiple time series}, Wiley, 1970 \\
\ni [HeLo] H. Helson and D. Lowdenslager,  Prediction theory and Fourier series in several variables, I, II.  {\sl Acta Math.} {\bf 99} (1958), 165-202, {\bf 106} (1961), 175-213.\\
\ni [HiPh] E. Hille and R. S. Phillips, {\sl Functional anaysis and semigroups}.  Colloq. Publ. 31, American Math. Soc., 1957. \\
\ni [In1] A. Inoue, Asymptotics for the partial autocorrelation function of a stationary process.  {\sl J. Analyse Math.} {\bf 81} (2000), 65-109. \\
\ni [In2] A. Inoue, Asymptotic behaviour for partial autocorrelation functions of fractional ARIMA processes.  {\sl Ann. Appl. Prob.} {\bf 12} (2002), 1471-1491. \\
\ni [In3] A. Inoue, AR and MA presentations of partial autocorrelation functions with applications.  {\sl Prob. Th. Rel. Fields} {\bf 140} (2008), 523-551. \\
\ni [InKa1] A. Inoue and Y. Kasahara, Partial autocorrelation functions of the fractional ARIMA process.  {\sl J. Multivariate Analysis} {\bf 89} (2004), 135-147. \\
\ni [InKa2] A.Inoue andY.Kasahara,\ Explicit\ representation\ of\ finite\ predictor\ coefficients  and its applications. {\sl Ann. Statist.} {\bf 34} (2006), 973-993. \\
\ni [JaLaEp] G. Janashia, E. Lagvilada and L. Ephremidze, A new method of matrix spectral factorization.  {\sl IEEE Trans. Info. Th.} {\bf 57} (2011), 2318-2326. \\
\ni [KaBi] Y. Kasahara and N. H. Bingham, Verblunsky coefficients and Nehari sequences.  Preprint, Hokkaido University. \\
\ni [LeMcK] N. Levinson and H. P. McKean, Weighted trigonometrical approximation on $R^1$ with application to the germ field of a stationary Gaussian noise.  {\sl Acta Math.} {\bf 112} (1964), 99-143 (reprinted in {\sl Selected papers of Norman Levinson} Vol. 2
(Birkh\"auser, Basel, 1998), IX, 222-266). \\
\ni [MaMiSc] A. Makagon, A. G. Miamee and B. S. W. Schr\"oder, Recursive condition for positivity of the angle for multivariate stationary sequences.  {\sl Proc. Amer. Math. Soc.} {\bf 126} (1998), 1821-1825. \\
\ni [Mas1] P. Masani, Crame\'er's theorem on monotone matrix functions and the Wold decomposition.  {\sl Probability and statistics: The Harald Cram\'er volume}
(ed. U. Grenander) 175-189, Wiley, 1959. \\
\ni [Mas2] P. Masani, The prediction theory of multivariate stochastic processes, III. {\sl Acta Math.} {\bf 104} (1960), 141-162. \\
\ni [Mas3] P. Masani, Shift-invariant spaces and prediction theory.  {\sl Acta Math.} {\bf 107} (1962), 275-290.\\
\ni [Mas4] P. Masani, Recent trends in multivariate prediction theory.  {\sl Multivariate Analysis} (Proc. Int. Symp., Dayton OH) 351-382, Academic Press, 1966.\\
\ni [Mas5] P. Masani, Comments on the prediction-theoretic ppers, 276-305 in [Wi]. \\
\ni [Mat] R. F. Matveev, Regularity of multi-dimensional stochastic processes with discrete time.  {\sl Dokl. Akad. Nauk SSSR} {\bf 126} (1959), 713-715. \\
\ni [MVT] M. Morf, A. Vieira and T. Kailath, Covariance characterization by partial autocorrelation matrices.  {\sl Ann. Statist.} {\bf 6} (1978), 643-678. \\
\ni [Nik] N. K. Nikolskii, {\sl Operators, functions and systems: an easy reading.  Volume 1:
Hardy, Hankel and Toeplitz; Volume 2: Model operators and systems}.  Math. Surveys and Monographs
{\bf 92, 93}, Amer. Math. Soc., 2002. \\
\ni [Pel1] V. V. Peller, Hankel operators and multivariate stochastic processes.  {\sl Proc. Symp. Pure Math.} 51, Part 1, 357-371, AMS, Providence, RI, 1990. \\
\ni [Pel2] V. V. Peller, Factorization and approximation problems for matrix functions.  {\sl J. American Math. soc.} {\bf 11} (1998), 751-770. \\
\ni [Pel3] V. V. Peller, {\sl  Hankel operators and their applications}, Springer, 2003.  \\
\ni [Pou1] M. Pourahmadi, A matricial extension of the Helson-Szeg\"o theorem and its application in multivariate prediction.  {\sl J. Multivariate Analysis} {\bf 16} (1985), 265-275. \\
\ni [Pou2] M. Pourahmadi, Joint mean-covariance models with applications to longitudinal data.  Unconstrained parametrization.  {\sl Biometrika} {\bf 86} (1999), 677-690.\\
\ni [Pou3] M. Pourahmadi, {\sl Foundations of time series analysis and prediction theory}.  Wiley, 2001.\\
\ni [Re] G. C. Reinsel, {\sl Elements of multivariate time series analysis}, Springer, 1997. \\
\ni [RiNa] F. Riesz and B. Sz.-Nagy, {\sl Le\c cons d'analyse fonctionnelle}, 2nd ed., Akad. Kiad\'o, 1953. \\
\ni [Ro1] Yu. A. Rozanov, Spectral properties of multivariate stationary processes and boundary properties of analytic functions.  {\sl Th. Probab. Appl.} {\bf 5} (1960), 362-376. \\
\ni [Ro2] Yu. A. Rozanov, {\sl Stationary random processes}, Holden Day, San Francisco, CA, 1967. \\
\ni [Sa1] D. Sarason, An addendum to "Past and future".  {\sl Math. Scand.} {\bf 30} (1972), 62-64.\\
\ni [Sa2] D. Sarason, Functions of vanishing mean oscillation.  {\sl Trans. Amer. Math. Soc.} {\bf 207} (1975), 391-405.\\
\ni [Sa3] D. Sarason, {\sl Function theory on the unit circle}.  Virginia Polytechnic Institute and State University, Blacksburg, VA, 1979.\\
\ni [SeSi] M. Seidel and B. Silbermann, Banach algebras of structured matrix sequences.  {\sl Linear Algebra Appl.} {\bf 430} (2009), 1243-1281. \\
\ni [Si1]   B. Simon, {\sl Orthogonal polynomials on the unit circle.  Part 1: Classical theory}.  AMS Colloq. Publ. 54.1, AMS, Providence, RI, 2005. \\
\ni [Si2]   B. Simon, {\sl Orthogonal polynomials on the unit circle.  Part 2: Spectral theory}.  AMS Colloq. Publ. 54.2, AMS, Providence, RI, 2005. \\
\ni [Si3] B. Simon, {\sl Szeg\"o's theorem and its descendants.  Spectral theory for} $L_2$ {\sl perturbatios of orthogonal polynomials}.  Princeton Univ. Press, Princeton, NJ, 2011. \\
\ni [Ste] E. M. Stein, {\sl Harmonic analysis}.  Princeton University Press. 1993. \\
\ni [Sto] M. H. Stone, {\sl Linear transformations in Hilbert space and their applications to analysis}.  Coll. Publ. XV. American Math. Soc., 1932. \\
\ni [SzNF] B. Sz.-Nagy and C. Foias, {\sl Harmonic analysis of operators on Hilbert space}, North-Holland, 1970
(2nd ed., with H. Bercovici and L. K\'erchy, Springer Universitext, 2010).\\
\ni [TrVo1] S. Treil and A. Volberg, Wavelets and the angle between past and future.  {\sl J. Functional Analysis} {\bf 143} (1997), 269-308.\\
\ni [TrVo2] S. Treil and A. Volberg, Completely regular multivariate stochastic processes and the Muckenhoupt conditiion.  {\sl Pacific J. Math.} {\bf 190} (1999), 361-382.\\
\ni [Wh] P. Whittle, On the fitting of multivariate autoregressions, and the approximate canonical factorization of a spectral density matrix.  {\sl Biometrika} {\bf 50} (1963), 129-134. \\
\ni [Wi] N. Wiener, {\sl Collected works, Volume III} (ed. P. Masani), MIT Press, 1981. \\
\ni [WiAk] N. Wiener and E. J. Akutowicz, A factorization of positive Hermitian matrices, {\sl J. Math. Mech.} {\sl 8} (1959), 111-120. \\
\ni [WiMa1] N. Wiener and P. Masani, The prediction theory of multivariate stochastic processes, I: The regularity condition. {\sl Acta Math.} {\bf 98} (1957), 111-150. \\
\ni [WiMa2] N. Wiener and P. Masani, The prediction theory of multivariate stochastic processes, II: The linear predictor. {\sl Acta Math.} {\bf 98} {\bf 99} (1958), 93-137. \\
\ni [WiMa3] N. Wiener and P. Masani, On bivariate stationary processes and the factorization of matrix-valued functions.  {\sl Th. Prob. Appl.} {\bf 4} (1959), 300-308.  \\
\ni [WiRo] R. A. Wiggins and E. A. Robinson, Recursive solution to the multichannel filtering problem.  {\sl J. Geophys. Res.} {\bf 70} (1965), 1885-1891.\\

\ni Mathematics Department, Imperial College, London SW7 2AZ, UK \\
\ni n.bingham@ic.ac.uk   nick.bingham@btinternet.com \\

\end{document}